# SL$(2,\mathbb{Z})$ – invariant spaces spanned by modular units

Wolfgang Eholzer and Nils-Peter Skoruppa

*Dedicated to the memory of Tsuneo Arakawa*

**Abstract**

Characters of rational vertex operator algebras (RVOAs) arising, *e.g.*, in 2 - dimensional conformal field theories often belong (after suitable normalization) to the (multiplicative) semi-group $E^+$ of modular units whose Fourier expansions are in $1+q\,\mathbb{Z}_{\geq 0}[\![q]\!]$, up to a fractional power of $q$. If, furthermore, all characters of a RVOA share this property then we have an example of what we call modular sets, *i.e.* finite subsets of $E^+$ whose elements (additively) span a vector space which is invariant under the usual action of SL$(2,\mathbb{Z})$. The appearance of modular sets is always linked to the appearance of other interesting phenomena. The first nontrivial example is provided by the functions appearing in the two classical Rogers-Ramanujan identities, and generalizations of these identities known from combinatorial theory yield further examples. The classification of modular sets and RVOAs seems to be related. This article is a first step towards the understanding of modular sets. We give an explicit description of the group of modular units generated by $E^+$, we prove a certain finiteness result for modular sets contained in a natural semi-subgroup $E_*$ of $E^+$, and we discuss consequences, in particular a method for effectively enumerating all modular sets in $E_*$.

## 1  Introduction

Two famous identities were discovered 1894 by Rogers [R] and rediscovered 1913 by Ramanujan and 1917 by Schur, and since then have been cited as

Rogers-Ramanujan identities:

$$\prod_{\substack{n \equiv \pm 1 \bmod 5 \\ n > 0}} (1 - q^n)^{-1} = \sum_{n \geq 0} \frac{q^{n^2}}{(1-q)(1-q^2)\cdots(1-q^n)},$$

$$\prod_{\substack{n \equiv \pm 2 \bmod 5 \\ n > 0}} (1 - q^n)^{-1} = \sum_{n \geq 0} \frac{q^{n^2+n}}{(1-q)(1-q^2)\cdots(1-q^n)}.$$

Apart from their combinatorial meaning concerning partitions, the Rogers-Ramanujan identities encode the following surprising fact. If we set $q = e^{2\pi i z}$ for $z$ in the complex upper half plane, and if we multiply the two identities by $e^{-\pi i z/30}$ and $e^{11\pi i z/30}$, respectively, then the functions involved in these identities become modular functions. As well-known and well-understood this statement appears for the products, which are, via the Jacobi triple product identities, quotients of elementary theta series, as remarkable this fact appears for the theta-like infinite series occurring in these identities. There is no known conceptual method in the theory of modular forms which produces modular functions of this shape[1].

The Rogers-Ramanujan identities are the first ones of an infinite series of identities of this kind, namely, of the Andrews-Gordon identities (see Section 2). The infinite series occurring in the Andrews-Gordon identities are more generally of the form

$$f_{A,b,c} = \sum_{n=(n_1,\ldots,n_r) \in \mathbb{Z}_{\geq 0}^r} \frac{q^{nAn^t + bn^t + c}}{(q)_{n_1} \cdots (q)_{n_r}}, \quad (q)_k = (1-q)(1-q^2)\cdots(1-q^k),$$

where $A$ is a symmetric positive definite rational $r \times r$-matrix, where $b$ is a rational row vector of length $r$ and $c$ is a rational number. Again, the $f_{A,b,c}$ occurring in the Andrews-Gordon identities are modular functions (since the products occurring in these identities are).

One may consider the following problem[2]: *For what $A$, $b$ and $c$ is $f_{A,b,c}$ a modular function ?*

As it turns out this seems to be a hard question: the answer is not known and the known instances of modular $f_{A,b,c}$ are very exceptional. For $r = 1$ the problem was completely solved by Zagier [Z]: there are precisely 7 triples of rational numbers $A > 0$, $b$ and $c$ such that $f_{A,b,c}$ is modular (see Table 1).

---
[1]At least no such method is known to the authors.
[2]The first one who mentioned this problem to the authors was Werner Nahm.

Table 1: All rational numbers $A > 0, b, c$ such that $f_{A,b,c}$ is modular.

| $A$ | 1 | 1 | $\frac{1}{2}$ | $\frac{1}{2}$ | $\frac{1}{2}$ | $\frac{1}{4}$ | $\frac{1}{4}$ |
|---|---|---|---|---|---|---|---|
| $b$ | 0 | 1 | 0 | $\frac{1}{2}$ | $-\frac{1}{2}$ | 0 | $\frac{1}{2}$ |
| $c$ | $-\frac{1}{60}$ | $\frac{11}{60}$ | $-\frac{1}{48}$ | $\frac{1}{24}$ | $\frac{1}{24}$ | $\frac{1}{40}$ | $\frac{1}{40}$ |

There are indications that a complete answer to this question would involve $K_3(\overline{Q})$ [N, Section 4], and might be related to the problem of classifying vertex operator algebras or two-dimensional quantum field theories [N, Section 4], [E-S].

However, the mentioned identities exhibit another remarkable fact. Namely, the space of linear combinations of the two products in the Rogers-Ramanujan identities is invariant under the natural action of $\mathrm{SL}(2,\mathbb{Z})$ on functions defined on the upper half plane. Moreover, the two products are modular units, they have non-negative integral Fourier coefficients and they are eigenfunctions under $z \mapsto z+1$. These properties hold also true for the products in the Andrews-Gordon identities (see section 2). More generally, such sets of products arise naturally as conformal characters of various (rational) vertex operator algebras [E-S]. The question of finding all such sets of modular units, the question about the modularity of the $f_{A,b,c}$, and the problem of classifying vertex operator algebras seem to be interwoven.

Hence, instead of trying to investigate directly the functions $f_{A,b,c}$ for modularity, one may hope to come closer to an answer to this problem by seeking first of all for a description of all finite sets of modular units of the indicated shape which span $\mathrm{SL}(2,\mathbb{Z})$-invariant spaces. We shall call such sets modular (see Section 2 for a precise definition). As it turns out, modular sets are indeed very exceptional and their description is a non-trivial task.

This article is first step towards the understanding of modular sets. As a byproduct we shall show that an important subclass of modular sets can be algorithmically enumerated.

## 2 Statement of results

A modular unit is a modular function on some congruence subgroup of $\Gamma := \mathrm{SL}(2,\mathbb{Z})$ which has no poles or zeros in the upper half plane $\mathbb{H}$.

Thus it takes on all its poles and zeros in the cusps. The set $U$ of all modular units is obviously a group with respect to the usual multiplication of modular functions.

In this note we are interested in modular units $f$ whose Fourier coefficients are non-negative integers, which satisfy $f(z+1) = c\,f(z)$ with a suitable constant $c$, and whose first Fourier coefficients are 1. Denote by $E^+$ the semi-group of all such units. In other words, $E^+$ is the semi-group of all modular units whose Fourier expansion is in $q^s(1+q\mathbb{Z}_{\geq 0}[\![q]\!])$ for some rational number $s$. Here $q^s$, for any real $s$, denotes the function $q^s(z) = \exp(2\pi i s z)$ with $z$ a variable in $\mathbb{H}$.

Special instances of $E^+$ are the units

$$[r]_l = q^{-l\mathbb{B}_2\left(\frac{r}{l}\right)/2} \prod_{\substack{n \equiv r \bmod l \\ n > 0}} (1-q^n)^{-1} \prod_{\substack{n \equiv -r \bmod l \\ n > 0}} (1-q^n)^{-1}, \qquad (1)$$

where $l \geq 1$ and $r$ are integers such that $l$ does not divide $r$ (*cf.* Lemma 6.1 in Section 6). Here we use $\mathbb{B}_2(x) = y^2 - y + \frac{1}{6}$ with $y = x - \lfloor x \rfloor$ as the fractional part of $x$.

In particular, we are interested in **modular sets**, by what we mean finite and non-empty subsets $S$ of $E^+$ such that the subspace (of the complex vector space of all functions on $\mathbb{H}$) which is spanned by the units in $S$ is invariant under $\Gamma$. Note that the group $U$ is invariant under $\Gamma$: if $f(z)$ is a unit and $A \in \Gamma$ then $f(Az)$ is again a unit. Thus it is easy to write down finite subsets of $U$ whose span is $\Gamma$-invariant. In contrast to this, $E^+$ is not invariant under $\Gamma$, and, indeed, as we shall explain in a moment, modular sets seem to be quite exceptional. We call a modular set nontrivial if it contains more units than merely the constant function 1.

An infinite series of examples for nontrivial modular sets is provided by the following. Let $l$ be an odd natural number and set

$$\phi_r = \prod_{\substack{1 \leq j \leq \frac{l-1}{2} \\ j \neq r}} [j]_l \qquad (1 \leq r \leq \tfrac{l-1}{2}).$$

Then, for each $l$, the set $\mathrm{AG}_l$ of all $\phi_r$ with $r$ in the given range is modular [E-S]. (See also [C-I-Z, Eq. (23)], where, however, the $\phi_r$ are not given as products, but as quotients of theta functions and the Dedekind $\eta$-function. Both expressions for the $\phi_r$ are easily identified on using the Jacobi triple product identity; *cf.* [E-S] for details.)

The existence of (nontrivial) modular sets is a somewhat remarkable fact. First of all, the notion of modular sets itself is bizarre: the action of $\Gamma$ on modular units defines automorphisms of the group of modular units, whereas a modular set requires the linear subspace, and not the subgroup, generated by its elements, to be $\Gamma$-invariant. More striking, modular sets seem to be bound to other remarkable phenomena. The functions $\phi_r$ occur as the one side of the Andrews-Gordon identities (see, *e.g.*, [B, Eq. (3.2), p. 15 ]):

$$\phi_r = q^{-l\mathbb{B}_2\left(\frac{r}{l}\right)/2} \sum_n \frac{q^{nAn^t + b_r n^t}}{(q)_{n_1} \cdots (q)_{n_{k-1}}}.$$

Here $k = (l-1)/2$ and $n = (n_1, \ldots, n_{k-1})$ runs over all vectors with non-negative integral entries, $A$ is the matrix $A = (\min(i, j))$, and $b_r$ is the vector with $\min(i+1-r, 0)$ as $i$-th entry, and finally

$$(q)_m = (1-q)(1-q^2) \cdots (1-q^m)$$

(with the convention $(q)_0 = 1$). The two identities for $l = 5$ are the classical Rogers-Ramanujan identities.

Finally, modular sets show up as sets of conformal characters of certain rational vertex operator algebras [E-S]. In fact, the modular sets $AG_l$ provide also examples of this [K-R-V, Eq. (2.1)-(2.3)], and we do not know any modular set which is not the set of conformal characters of a rational vertex operator algebra [E-S].

The ultimate goal would be a classification of all modular sets. As indicated in [E-S] this is related to the open problem of the classification of a certain class of rational vertex operator algebras arising in 2-dimensional conformal field theories.

The first natural step in the study of modular sets is to ask for a more explicit description of the semi-group $E^+$. We shall prove the following structure theorem.

**Theorem 1.** *Let $E$ be the group of units generated by the $[r]_l$ (defined in (1)). Then $\mathbb{Q}^* \cdot E$ coincides with the group of all modular units whose Fourier expansions are in $q^s \mathbb{Q}[\![q]\!]$ for suitable rational numbers $s$.*

In particular, the group $E$ is identical with the group of modular units whose Fourier expansion is in $q^s(1 + q\mathbb{Z}[\![q]\!])$ for some rational number $s$. Thus, the group of modular units generated by $E^+$ is obviously contained in $E$. Since it contains on the other side the generators $[r]_l$ of $E$, we conclude

**Corollary to Theorem 1.** *The group of modular units generated by the elements of $E^+$ coincides with the one generated by the $[r]_l$. In particular, each element of $E^+$ is a product of integral, though not necessarily positive, powers of the special units $[r]_l$.*

There is another remarkable consequence of Theorem 1. Namely, the first Fourier coefficient of a conformal character needs not to be 1. Thus, with regard to applications to conformal characters, it would be more natural to study the semi-group of modular units with Fourier expansions in $q^s \, \mathbb{Z}_{\geq 0}[\![q]\!]$ for some $s$. However, by the theorem this semi-group equals $\mathbb{Z}_{>0} \cdot E^+$, which shows that one does not loose any generality by restricting to $E^+$, as we *a priori* did in this article.

It is worthwhile to describe the structure of the group $E$, *i.e.* the (multiplicative) relations satisfied by the generators $[r]_l$ of $E$. For each $l$ we have the obvious homomorphism $\mathbb{Z}[\mathbb{Z}/l\mathbb{Z}] \to E$, which associates to a $\mathbb{Z}$-valued map $f$ on $\mathbb{Z}/l\mathbb{Z}$ the product of all $[r]_l^{f(r)}$, where $r$ runs through a complete set of representatives for the nonzero residue classes modulo $l$. Moreover, one easily verifies the the distribution relations

$$[r]_l = \prod_{\substack{s \bmod m \\ s \equiv r \bmod l}} [s]_m,$$

valid for all $l$ and $m$ such that $l|m$. We may thus combine the above homomorphisms by setting, for any locally constant $f \colon \widehat{\mathbb{Z}} \to \mathbb{Z}$ with $f(0) = 0$,

$$[\ ]^f := \prod_{r \bmod l} [r]_l^{f(r)}.$$

Here $\widehat{\mathbb{Z}}$ denotes the Pruefer ring, (*i.e.* $\widehat{\mathbb{Z}} = \text{proj}\lim \mathbb{Z}/l\mathbb{Z}$, equipped with the topology generated by the cosets $\widehat{\mathbb{Z}}/l\widehat{\mathbb{Z}}$), and $l$ is any positive integer such that $f$ is constant on the cosets modulo $l\widehat{\mathbb{Z}}$. By the distribution relations $[\ ]^f$ does not depend on a particular choice of $l$. One has:

**Supplement to Theorem 1.** *The map $f \mapsto [\ ]^f$ induces an isomorphism of $\mathrm{L}(\widehat{\mathbb{Z}})/\mathrm{L}(\widehat{\mathbb{Z}})^-$ and the group $E$ of modular units generated by the $[r]_l$ (defined in (1)). Here $\mathrm{L}(\widehat{\mathbb{Z}})$ is the group of $\mathbb{Z}$-valued, locally constant maps on $\widehat{\mathbb{Z}}$ vanishing at 0, and $\mathrm{L}(\widehat{\mathbb{Z}})^-$ is the subgroup of odd maps.*

The division by $\mathrm{L}(\widehat{\mathbb{Z}})^-$ is due to the (obvious) relations $[r]_l = [-r]_l$.

Denote by $E_*$ the semi-group of products of non-negative powers of the special functions $[r]_l$. Clearly, $E^+$ contains the semi-subgroup $E_*$, and, by the corollary to Theorem 1, $E^+$ and $E_*$ generate the same group. However, $E^+$ is strictly larger than $E_*$; e.g., the function $[1]_4^3/[2]_4 = \eta^{-1} \sum_n q^{n^2}$ (with $\eta$ denoting the Dedekind eta-function) is in $E^+$, but not in $E_*$. Understanding the last example and giving a complete description of $E^+$ seems to be difficult.

Therefore, we shall consider in the following only modular subsets which are contained in the the semi-subgroup $E_*$ of products of non-negative powers of the $[r]_l$. This restriction seems to be not too serious: in fact, the only examples of modular sets not contained in $\mathbb{Z}_{>0} \cdot E_*$ which we know are in a certain sense trivial (*cf.* [E-S]).

As the second main result of the present article, we shall prove a certain finiteness property for modular subsets of $E_*$, which will in particular imply a method to systematically enumerate them. Namely, for fixed positive integers $n$ and $l$, let $E_n(l)$ be the set of all products of the form

$$[r_1, \ldots, r_k]_l := \prod_{j=1}^{k} [r_j]_l,$$

with $k \leq n$, and arbitrary integers $r_j$ which are not divisible by $l$. The sets $E_n(l)$ are clearly finite. Using the distribution relations it is clear that *any modular subset of $E_*$ is contained in some $E_n(l)$ with suitable $n$ and $l$*. We shall prove:

**Theorem 2.** *For each $n$ the number of $l$ such that $E_n(l)$ contains a nontrivial modular set is finite. More precisely, if $E_n(l)$ contains a nontrivial modular set, then $l \leq 13.7^n$.*[3]

Our proof will exhibit a method to compute, for a given $n$, all modular subsets of $E_n(l)$ for all $l$. This method, however, becomes quickly non-realistic for growing $n$.

In Table 2 we listed all modular subsets of $E_n(l)$ for $n \leq 3$ and $l \geq 1$. For each $n$, we listed only those modular sets which do not already belong to some $E_k(l)$ with $k < n$, and which cannot be decomposed into a disjoint union of smaller modular sets. By $S^n$, for a modular set $S$ and a positive integer $n$, we denote the set of all $n$-fold products of functions in $S$. Obviously, $S^n$ is again modular. Note that, for $n \leq 3$, there is exactly one 'new' modular set, which we called $W_7$. More examples of modular sets can be found in [E-S].

---

[3] Actually the existence of a nontrivial modular set in $E_n(l)$ implies $l \leq 5n$. However, this sharper result relies on a deep analysis of the (projective) $\mathrm{SL}(2,\mathbb{Z})$-module of all modular forms of weight $\frac{1}{2}$, and will be published elsewhere.

Table 2: All modular subsets of $E_n(l)$ for $n \leq 3$ and arbitrary $l$.

| | $l = 5$ | 7 | 9 |
|---|---|---|---|
| $n = 1$ | $AG_5$ | | |
| 2 | $AG_5^2$ | $AG_7$ | |
| 3 | $AG_5^3$ | $W_7 := \{[1,2,3]_7\} \cup \{[r,r,3r]_7 : r = 1,2,3\}$ | $AG_9$ |

The plan of the rest of this article is as follows: In Section 3 we shall prove Theorem 1 and its supplement, and in Section 4 we shall prove Theorem 2. The auxiliary results derived in Section 4 have some interest, independent of the proof of Theorem 2, in connection with the question of searching for modular sets. In Section 5 we shall briefly indicate how to use these auxiliary results for calculating, *e.g.*, the above table.

In the proofs of the two theorems we need certain properties of the $[r]_l$'s, which we derive in Section 6 by rewriting $[r]_l$ in terms of $l$-division values of the Weierstrass $\sigma$-function and using some of their basic properties. Since we did not find any convenient reference to cite these properties directly we decided to develop quickly from scratch the corresponding theory in form of a short Appendix and part of Section 6. In particular, we emphasize in the Appendix that the Weierstrass $\sigma$-function and its $l$-division values are best understood by viewing the Weierstrass $\sigma$-function as a Jacobi form on the full modular group of weight and index equal to $\frac{1}{2}$ (see Theorem 7.1).

## 3 The group of units generated by the $[r]_l$

In this section we prove Theorem 1 and its supplement. We shall actually prove the slightly stronger Theorem 3.2. Its proof depends on two well-known facts: first, that the group of all modular units modulo the so-called Siegel units is a torsion group, and, secondly, that modular forms on congruence subgroups with rational Fourier coefficients have bounded denominators. The short proof of the first one is given in Section 6, for the second, deeper one, we refer to the literature.

We precede the proof of Theorem 3.2 by three lemmas. The first one,

which we actually call theorem to emphasize its more general usefulness, is a general statement about product expansions of holomorphic and periodic functions in the upper half plane. It is important for the proof of the third lemma, but it also implies directly the supplement to Theorem 1 of Section 2.

**Theorem 3.1.** *Let $f$ be a holomorphic and periodic function on the upper half plane whose Fourier expansion is in $1 + q\mathbb{Z}[\![q]\!]$ . Then there exists a unique sequence $\{a(n)\}$ of integers such that*

$$f = \prod_{n \geq 1}(1 - q^n)^{a(n)},$$

*for sufficiently small $|q|$.*

*Remark* 3.1. As can be read off from the proof the lemma actually holds true with $\mathbb{Z}$ replaced by an arbitrary subring of $\mathbb{C}$.

*Proof.* The existence of the sequence $a(n)$ follows by induction on $n$. Namely, assume that one has already found integers $a(n)$ $(1 \leq n < N)$ such that

$$g := f / \prod_{n=1}^{N-1}(1 - q^n)^{a(n)} = 1 + \mathcal{O}(q^N).$$

Let $-a(N)$ be the Fourier coefficient of $g$ in front of $q^N$. Clearly $a(N)$ is integral. One has

$$f / \prod_{n=1}^{N}(1 - q^n)^{a(n)} = g/(1 - q^N)^{a(N)} = 1 + \mathcal{O}(q^{N+1}).$$

The uniqueness of the $a(n)$ follows from the uniqueness of the Fourier expansion of $q\frac{d}{dq} \log f$. □

*Proof of Supplement to Theorem 1.* That the kernel of the map $\mathrm{L}(\widehat{\mathbb{Z}}) \mapsto E$ equals $\mathrm{L}(\widehat{\mathbb{Z}})^-$ follows from the uniqueness of the product expansion in the preceding proposition and on writing

$$[\,]^f = \prod_{r \bmod l} [r]_l^{f(r)} = q^c \prod_{n \geq 1}(1 - q^n)^{-f(n)-f(-n)}$$

with a suitable constant $c$. The surjectivity is obvious from the definition of $E$. □

**Lemma 3.1.** *Let $f \in \frac{1}{D}\mathbb{Z}[\![q]\!]$ for some positive integer $D$. If some positive integral power of $f$ has integral coefficients, then $f$ has integral coefficients.*

*Proof.* By assumption about the coefficients of $f$ we can write $f = \gamma \cdot h$ with a suitable rational number $\gamma$ and with a primitive $h$. Here primitive means that $h$ is a power series in $q$ with integral coefficients $a(l)$ which are relatively prime. By assumption, $\gamma^N \cdot h^N$, for some integer $N \geq 1$, has integral coefficients. We shall show in a moment that $h^N$ is primitive. From this we deduce that $\gamma^N$ is integral. Hence $\gamma$ is integral, which proves the lemma.

It remains to show that $h^N$ is primitive. Let $p$ be a prime. Since $h$ is primitive, there exists an $l$ such that $p|a(j)$ for $j < l$ and $p \nmid a(l)$. But then the $q^{Nl}$-coefficient of $h^N$ satisfies

$$\sum_{i_1+\cdots+i_N=Nl} a(i_1)\cdots a(i_N) \equiv a(l)^N \bmod p,$$

and whence is not divisible by $p$. $\square$

For the following, let $E(l)$, for fixed $l$, denote the group generated by the $[r]_l$ with $1 \leq r \leq \lfloor l/2 \rfloor$.

**Lemma 3.2.** *Let $f$ be a modular unit with rational Fourier coefficients. Assume that a positive integral power of $f$ lies in $E(l)$. Then $f$ is in $E(2l)$ (and even in $E(l)$ for odd $l$).*

*Proof.* Since $f$ is invariant under a congruence subgroup it has bounded denominators, *i.e.* there exist an integer $D > 0$ such that $D \cdot f$ has integral Fourier coefficients. This well-known fact follows, *e.g.*, on writing $f\eta^{24N}$, with a suitable integer $N > 0$ (and with $\eta$ denoting the Dedekind eta-function), as linear combination of modular forms with integral Fourier coefficients (which is possible by Theorem 3.52 in [Sh]), deducing from this that $f\eta^{24N}$ has bounded denominators, which in turn implies that $f$ has bounded denominators since $\eta^{-1}$ has integral Fourier coefficients.

Combining the latter with the fact that some positive integral power of $f$ lies in $E(l)$, we see that, for some rational number $s$, the function $q^{-s}f$ satisfies the assumption of Lemma 3.1, and hence is in $\mathbb{Z}[\![q]\!]$. Moreover, by assumption, its first Fourier coefficient is 1.

But then $q^{-s}f$ possesses a product expansion as in the Theorem 3.1. By the uniqueness of the $a(n)$, and since a nonzero integral power $f^N$ of $f$ is a product of $[r]_l$'s, we conclude that $Na(n) = Na(m)$ for $n \equiv \pm m \bmod l$, and that $Na(0) = 0$. Since $N \neq 0$ the same holds true with $\{Na(n)\}$ replaced by $\{a(n)\}$. Thus we find, on re-ordering the product expansion of $f$ according to the residue classes of $n$ modulo $l$, that

$$f^{-1} = [1]_l^{a(1)}[2]_l^{a(2)} \cdots [m-1]_l^{a(m-1)}[m]_l^{\nu a(m)}$$

where $m = \lfloor l/2 \rfloor$, and where $\nu = 1$ for odd $l$, and $\nu = 1/2$ for even $l$. Hence, if $l$ is odd, then $f \in E(l)$. If $l$ is even, one uses the distribution relations (in particular, $[l/2]_l^{1/2} = [l/2]_{2l}$) to deduce $f \in E(2l)$. □

**Theorem 3.2.** *The group of modular units on $\Gamma(l)$ $(= \{A \in \mathrm{SL}(2, \mathbb{Z}) \mid A \equiv 1 \bmod l\})$ with Fourier expansions in $q^s \mathbb{Q}[\![q]\!]$ for suitable rational numbers $s$ is a subgroup of $\mathbb{Q}^* \cdot E(2l)$ (and even of $\mathbb{Q}^* \cdot E(l)$, for odd $l$).*

*Proof.* Let $f$ be unit on $\Gamma(l)$ such that, for some rational number $s$, the function $q^{-s}f$ has rational Fourier coefficients. For showing that $f$ is contained in $\mathbb{Q}^* \cdot E(l)$ or $\mathbb{Q}^* \cdot E(2l)$, respectively, we may assume that $f$ is normalized, i.e. that its first Fourier coefficient is 1. By Lemma 3.2 it then suffices to show that a positive integral power of $f$ lies in $E(l)$.

By Theorem 6.2 of Section 6 we know that some nontrivial power of $f$ can be written as product of Siegel units $s_\alpha$, which are defined by Eq. (2) of Section 6. More precisely, there exists integers $a > 0$, $b(\alpha)$, and a constant $c$ such that

$$f^a = c \cdot \prod_{\alpha \in I} s_\alpha^{b(\alpha)},$$

where $I$ is a finite set of pairs of rational numbers of the form $(\frac{r}{l}, \frac{s}{l})$ with integers $r, s$ such that $\gcd(r, s, l) = 1$.

By replacing $a$ and the $b(\alpha)$ by suitable positive integral multiples we may assume that $f^a$ is invariant under $T = (1, 1; 1, 0)$. On the other hand, by Theorem 6.1 in Section 6 we have that $s_\alpha \circ T$ equals $s_{\alpha T}$, up to multiplication by a constant. Let $K$ denote the field of $l$-th roots of unity. For an integer $y$ relatively prime to $l$ denote by $\sigma_y$ the automorphism of $K$ which maps an $l$-th root of unity $\zeta$ to $\zeta^y$. We extend $\sigma_y$ to an automorphism of the ring $R = \bigoplus_{s \in \mathbb{Q}} q^s K[\![q]\!]$ by letting it act on coefficients. Since $f$ has rational coefficients, it is invariant under $\sigma_y$. From the formula for $s_\alpha$ in Section 6 it is immediate that, for $\alpha \in I$, one has $s_\alpha \in R$ and that $\sigma_y s_\alpha$ equals $s_{\alpha D(y)}$, up to multiplication by a constant and with $D(y) = (1, 0; 0, y)$.

Using these properties we can write

$$f^{al\varphi(l)} = \prod_{y \bmod *l} \prod_{h=0}^{l-1} \sigma_y(f^a \circ T^h) = d \prod_{\alpha \in I} \Big( \prod_{y \bmod *l} \prod_{h=0}^{l-1} s_{\alpha T^h D(y)} \Big)$$

with a suitable constant $d$, where the asterisk indicates that $y$ runs through a complete set of primitive residue classes modulo $l$, and $\varphi(l)$ denotes as usual the number of such classes.

It remains to show that the expressions $t_\alpha$ in the rightmost parenthesis are in $E(l)$, up to multiplication by constants (whose product then equals $d^{-1}$, since $f$ is normalized). Write $\alpha = (r,s)/l$ as above. Clearly $\alpha T^h D(y) = (r,t)/l$ with a suitable integer $t$. If $h$ and $y$ run through the given range, then t runs through a complete set of representatives for the residue classes modulo $l$ which are relatively prime to $\gcd(r,l)$, and each such $t$ is taken on the same number of times, say $p$ (look at the action of the subgroup of $\mathrm{GL}(2,\mathbb{Z}/l\mathbb{Z})$ of matrices of the form $(1,x;0,y)$ on pairs of residue classes $(u,v)$ in $(\mathbb{Z}/l\mathbb{Z})^2$ with $\gcd(u,v,l)=1$ ).

Thus $t_\alpha$ is the $p$-th power of

$$\prod_{\substack{t \bmod l \\ \gcd(t,r,l)=1}} s_{(r,t)/l} = \prod_{t \bmod l} \prod_{d|r,l,t} s_{(r,t)/l}^{\mu(d)} = \prod_{d|r,l} \prod_{u \bmod l/d} s_{(\frac{r}{d},u)/\frac{l}{d}}^{\mu(d)} = \prod_{d|r,l} \left[\frac{r}{d}\right]_{\frac{l}{d}}^{\mu(d)}.$$

Here we used the Moebius function $\mu(d)$, and, for the last identity, Lemma 6.1 of Section 6; moreover, we have to assume that $l$ does not divide $r$ (since $s_\alpha$, for $\alpha \in \mathbb{Z}^2$, is not defined). On using the distribution relations in $E$ we can rewrite the right hand side as power products of $[r]_l$'s. If $l$ divides $r$, then we leave it to the reader to verify by a similar calculation (using directly the definition (2) of $s_\alpha$) that the left hand side of the last identity equals $\prod_r [r]_l$, where $r$ runs through a complete system of representatives for the primitive residue classes modulo $l$. □

*Proof of Theorem 1.* This is clearly a consequence of Theorem 3.2. □

## 4  Properties of modular sets

In this section we shall prove Theorem 2. Actually, we shall prove the slightly stronger Theorem 4.2 below. Its proof will mainly depend on two results:

the first one concerns a sort of measure on the projective space over $\mathbb{Z}/l\mathbb{Z}$ (Theorem 4.1; see also the beginning of §5). The second result (Lemma 4.2) uses information about the action of $\Gamma$ on the $[r]_l$, and will not be completely proved before Section 6.

The first lemma gives a necessary criterion for a set $S \subset E^+$ to be modular in terms of the vanishing or pole orders of the functions in $S$. Let $f \not\equiv 0$ be a modular function on some subgroup of $\Gamma$, and let $s \in \mathbb{P}^1(\mathbb{Q}) = \mathbb{Q} \cup \{\infty\}$ be any cusp. Then there exists a $A \in \Gamma$ such that $s = A\infty$, and a real number $\alpha$ such that $f(Az)q^{-\alpha}(z)$ tends to a non-zero constant for $z = it$ with real $t \to \infty$. The number $\alpha$ does not depend on the choice of $A$. We set

$$\mathrm{ord}_s(f) = \alpha.$$

**Lemma 4.1.** *Let $S$ be a finite set of modular functions such that the space spanned by its elements is invariant under $\mathrm{SL}(2,\mathbb{Z})$. Then the map*

$$\nu : \mathbb{P}^1(\mathbb{Q}) \to \mathbb{Q}, \quad \nu(s) = \min_{f \in S} \mathrm{ord}_s(f)$$

*is constant.*

*Proof.* Indeed, for any fixed $A, B \in \mathrm{SL}(2,\mathbb{Z})$ and any $f \in S$ the function $f \circ A$ is a linear combination of the functions $g \circ B$ with $g \in S$. In particular, comparing the leading terms of the Fourier expansions of these functions, we conclude

$$\mathrm{ord}_\infty(f \circ A) \geq \min_{g \in S} \mathrm{ord}_\infty(g \circ B).$$

Since this is true for any $f$, and on using $\mathrm{ord}_\infty(f \circ A) = \mathrm{ord}_{A\infty}(f)$ we obtain $\nu(A\infty) \geq \nu(B\infty)$. Interchanging the role of $A$ and $B$ we see that here we actually have an equality. This proves the lemma. $\square$

**Lemma 4.2.** *Let $s \in \mathbb{P}^1(\mathbb{Q})$. Then*

$$\mathrm{ord}_s([r]_l) = -\frac{t^2}{2l} B_2\left(\frac{ar}{t}\right),$$

*where $s = \frac{a}{c}$ with relatively prime integers $a$ and $c$ (in particular, $a = \pm 1$ and $c = 0$, if $s = \infty$), and where $t = \gcd(c, l)$.*

*Proof.* Let $A \in \Gamma$ be a matrix with first row equal to $(a, c)^t$, i.e. such that $A\infty = s$. Then $\mathrm{ord}_s([r]_l) = \mathrm{ord}_\infty([r]_l \circ A)$, and the right hand side is given in Lemma 6.1 of Section 6. $\square$

Combining Lemma 4.1 and Lemma 4.2 we obtain the following necessary criterion for a set $S \subset E_n(l)$ to be modular. This criterion is the key for the proof of Theorem 2. We remark that Lemma 4.3 is actually the only instance in the proof of Theorem 2, where we use that $S$ is contained in $E_*$, rather than only in $E^+$

**Lemma 4.3.** *Let $S \subset E_n(l)$ be modular, and assume that $S$ contains at least one $n$-fold product (i.e. an element in $E_n(l) \setminus E_{n-1}(l)$). Then, for all divisors $t$ of $l$, one has*

$$\max_{[a_1,\ldots,a_k]_l \in S} \max_{a \bmod {}^*t} \sum_{j=1}^{k} \mathbb{B}_2\left(\frac{a_j a}{t}\right) = \frac{n}{6t^2}.$$

*Here the asterisk indicates that $a$ runs through a complete set of representatives for the primitive residue classes modulo $t$.*

*Remark* 4.1. Note that the lemma implies that $\gcd(a_1,\ldots,a_n,l) = 1$ for all $n$-fold products $\pi = [a_1,\ldots,a_n]_l \in S$. Indeed, if $d$ denotes this gcd, then the lemma applied to $t = d$ becomes $\frac{n}{6} \leq \frac{n}{6d^2}$, whence $d = 1$.

*Proof.* If $f \in S$ is an $k$-fold product, then $\mathrm{ord}_0(f) = -\frac{k}{12l}$ by the preceding lemma. Since $S$ contains an $n$-fold product, we conclude that

$$\min_{f \in S} \mathrm{ord}_0(f) = -\frac{n}{12l}.$$

The claimed inequality is now an immediate consequence of the first two lemmas. □

We call a point $P = (\bar{a}_1,\ldots,\bar{a}_n) \in (\mathbb{Z}/l\mathbb{Z})^n$ *special* if

$$\sum_{j=1}^{n} \mathbb{B}_2\left(\frac{a_j a}{t}\right) \leq \frac{n}{6t^2}$$

for all divisors $t$ of $l$ and all integers $a$ relatively prime to $t$. Here the bar denotes reduction modulo $l$.

Theorem 2 will now be a consequence of Lemma 4.3 and the following theorem, whose proof will take the rest of this section.

**Theorem 4.1.** *For a given $n$ there exist only a finite number of $l$ such that $(\mathbb{Z}/l\mathbb{Z})^n$ contains a special point. More precisely, if $(\mathbb{Z}/l\mathbb{Z})^n$ with $l > 1$ contains a special point, then*

$$l \leq B := \left( \frac{2\left(1 + l^{\frac{1}{n}-1}\right)}{1 - \sqrt{\frac{1}{3} + \frac{2}{3p^2}}} \right)^n ,$$

*where $p$ is the smallest prime divisor of $l$.*

**Theorem 4.2.** *If $E_n(l)$ contains a nontrivial modular set, then $l \leq B$ with $B$ as in Theorem 4.1.*

*Proof.* Let $S$ be a modular subset of $E_n(l)$, and let $k$ be minimal such that $S$ is contained in $E_k(l)$. Let $\pi = [r_1, \ldots, r_k] \in S$. By Lemma 4.3 $\pi$ yields a special point $(\overline{r}_1, \ldots, \overline{r}_k) \in (\mathbb{Z}/l\mathbb{Z})^k$. Hence, by Theorem 4.1, $l$ is bounded from above by the right hand side of the claimed inequality, but with $n$ replaced by $k$. Since $k \leq n$, the theorem then follows. □

*Proof of Theorem 2.* This is an immediate consequence of the preceding theorem. The bound in Theorem 2 is obtained from the bound of Theorem 4.2 by estimating $p$ to below by 2 and on using $1 + l^{1/n-1} \leq 2$. □

It remains to prove Theorem 4.1 on special points. For its proof we use

**Lemma 4.4.** *Let $P \in (\mathbb{Z}/l\mathbb{Z})^n$. Then there exists an integer $b$ not divisible by $l$ such that $b \cdot P = (\overline{b}_1, \ldots, \overline{b}_n)$ with integers $b_j$ (and where the bar denotes reduction modulo $l$) satisfying*

$$|b_j| \leq l^{1-\frac{1}{n}} + 1.$$

*Remark* 4.2. Note that the inequality is, for fixed $n$ and asymptotically in growing primes $l$, best possible, apart from a constant. Indeed, the number of points in $(\mathbb{Z}/l\mathbb{Z})^n$ described by homogeneous coordinates satisfying the above inequality is

$$\leq (2l/l^{1/n} + 3)^n \approx 2^n l^{n-1}.$$

But, for growing primes $l$, this is up to factor $2^n$ asymptotically equal to the number of orbits of $(\mathbb{Z}/l\mathbb{Z})^n$ modulo multiplication by non-zero elements of $\mathbb{Z}/l\mathbb{Z}$, which is

$$\frac{(l^n - 1)}{l - 1} + 1.$$

*Proof.* For an integer $r$, set $B_r = [-r, r]^n \cap \mathbb{Z}^n$, and let $C_r$ denote the reduction of $B_r$ modulo $l$. Assume $r < \frac{l}{2}$. Then $C_r$ contains exactly $(2r+1)^n$ elements. Note that the sum of two points of $C_r$ always lies in $C_{2r}$.

Consider the sets $x \cdot P + C_r$, where $x$ runs through $\mathbb{Z}/l\mathbb{Z}$. If the sum of the cardinalities of these sets is strictly greater than $l^n$, i.e. if $l \cdot (2r+1)^n > l^n$, then there exist at least two which have non-empty intersection.

Assume that there exists an integer $r$ satisfying the inequalities of the two preceding paragraphs, *i.e.* satisfying

$$\frac{l}{2} > r > \frac{l^{1-1/n}}{2} - \frac{1}{2} =: \rho.$$

Pick $x \not\equiv x' \bmod l$ such that $x \cdot P + C_r$ and $x' \cdot P + C_r$ contain a common point $Q$. Then $xP - Q$ and $Q - x'P$ both lie in $C_r$, and hence their sum $(x - x')P$ is in $C_{2r}$, whence can be represented by a point in $B_{2r}$.

If $\frac{l}{2} > \rho + 1$ we may take $r = \lfloor \rho + 1 \rfloor$ to fulfill the above two inequalities. Since then $2r \leq 2\rho + 2 = l^{1-1/n} + 1$, the lemma follows. Otherwise $\frac{l}{2} \leq \rho + 1 \leq 2\rho + 2 = l^{1-1/n} + 1$, and then the lemma is trivial. $\square$

*Proof of Theorem 4.1.* Let $P \in (\mathbb{Z}/l\mathbb{Z})^n$. Choose $b$ as in the last lemma. Write $\frac{b}{l} = \frac{a}{t}$ with a divisor $t$ of $l$ and $\gcd(a, t) = 1$. Note that $t \neq 1$ (since $b$ is not divisible by $l$), and hence $t \geq p$ with the smallest prime divisor $p$ of $l$. Thus,

$$\frac{aP}{t} \equiv \frac{1}{l}(b_1, \ldots, b_n) \bmod \mathbb{Z}^n,$$

with integers $b_j$ satisfying

$$|b_l/l| \leq l^{-\frac{1}{n}} + l^{-1} =: s.$$

Since $\mathbb{B}_2$ is decreasing in $[0, \frac{1}{2}]$, we find, for $s \leq \frac{1}{2}$, *i.e.* for $l \geq (2(1 + l^{\frac{1}{n}-1})^n$, the inequality

$$\sum_{j=1}^{n} \mathbb{B}_2\left(\frac{ab_j}{t}\right) \geq n\mathbb{B}_2(s).$$

Thus, if $s$ satisfies

$$\mathbb{B}_2(s) > \frac{1}{6p^2} (\geq \frac{1}{6t^2}),$$

then $(\mathbb{Z}/l\mathbb{Z})^n$ can never contain a special point. It is easily checked that the last inequality, together with $l^{-\frac{1}{n}} + l^{-1} = s \leq \frac{1}{2}$, is equivalent to

$$s < \frac{1 - \sqrt{\frac{1}{3} + \frac{2}{3p^2}}}{2}.$$

From this the theorem becomes obvious. □

## 5 Computing modular sets

We explain how we computed Table 2 in Section 2. The remarks of this section can actually be used to find for arbitrary $n$ all modular sets in $E_n(l)$ for all $l$, though, for growing $n$, the required computational resources become soon unrealistic[4].

The algorithm to enumerate all modular sets is based on the simple observation, that a subset of $E_n(l)$ is modular if and only if the space of functions which it spans is invariant under $z \mapsto -1/z$. This characterization follows easily from the fact that $\mathrm{SL}(2,\mathbb{Z})/\{\pm 1\}$ is generated by $z \mapsto -1/z$ and $z \mapsto z+1$. The invariance of a space spanned by a subset of $E_n(l)$ can be checked using explicit transformation formulas for the $[r]_l$ under $\Gamma$; *cf.* Theorem 6.1 and Lemma 6.1. (Note that by standard arguments from the theory of modular forms it suffices to check $N = N(l,n)$ many Fourier coefficients only to decide whether $\pi(-1/z)$, for $\pi$ in $E_n(l)$, is a linear combination of products in $E_n(l)$, where $N(l,n)$ is a constant depending only on $l$ and $n$, and which can be determined explicitly.)

Now, to find the maximal modular subset of $E_n(l)$ one could proceed as follows: Compute the set $S_1$ of all $\pi$ in $E_n(l)$ such that $\pi(-1/z)$ is a linear combination of functions in $E_n(l)$. Next one computes the set $S_2$ of all functions $\pi$ in $S_1$ such that $\pi(-1/z)$ is a linear combination of functions in $S_1$. Continuing like this one obtains a decreasing sequence of sets $S_k$. Either at some point $S_k$ is empty, and then $E_n(l)$ contains no modular subset, or else $S_k = S_{k+1}$ for some $k$, and then $S_k$ is the maximal modular subset of $E_n(l)$. However, since the number of products in $E_n(l)$ grows exponentially in $l$ it is necessary to look theoretical means to reduce the computational complexity. We indicate two such means.

---

[4]Using a cluster all modular sets in $E_n(l)$ for all $l$ and $n \leq 13$ could be determined; the results of this computation will be published elsewhere.

Assume $S \subset E_n(l)$ is modular. In this paper we are only interested in $n < 6$. This simplifies the computations a bit since then, for each $0 \leq k \leq n$, the subset $S(k)$ of all products of length $k$ in $S$ is already modular. Indeed, by the very definition of $[r]_l$ the elements of $S(k)$ have a Fourier expansion in powers of $q^{n/12l}$, where $n \in -kl^2 + 6\mathbb{Z}$. Furthermore, the $\pi(-1/z)$, for $\pi \in S(k)$, have a Fourier expansion in powers of $q^{n/12l}$, where $n \in -k + 12\mathbb{Z}$ (*cf.* Lemma 6.1). From this our proposition follows immediately if 3 does not divide $l$ (and since $k < 6$). If 3 divides $l$, then our argument shows that $S$ can only contain products of length 3 (again using $k < 6$), and our claim follows in this case too.

Hence, for verifying our table we can restrict our search for modular subsets of $E_n(l)$ to modular subsets of $F_n(l) := E_n(l) \setminus E_{n-1}(l)$.

Next, it is not at all necessary to consider all functions in $F_n$. Namely, let us call a subset $T$ of $\mathbb{P}^{n-1}(\mathbb{Z}/l\mathbb{Z})$ premodular if
$$\max_{P \in T} \beta_t(P) = \frac{n}{6t^2}$$
for all divisors $t$ of $l$. Here we use
$$\beta_t([\overline{a}_1 : \cdots : \overline{a}_n]) = \max_{a \bmod {}^* t} \sum_{j=1}^n \mathbb{B}_2\left(\frac{aa_j}{t}\right)$$
(with the asterisk as in Lemma 4.3 and the bar denoting reduction modulo $l$). Let $C_n(l)$ be the union of all premodular subsets in $\mathbb{P}^{n-1}(\mathbb{Z}/l\mathbb{Z})$, if there are any, and $C_n(l) = \emptyset$ otherwise.

If $S \subset F_n(l)$ is modular, then, by Lemma 4.3, the set $\overline{S}$ of all points $[\overline{a}_1 : \cdots : \overline{a}_n] \in \mathbb{P}^{n-1}(\mathbb{Z}/l\mathbb{Z})$ such that $\pi = [a_1, \ldots, a_n]_l \in S$ is premodular.

Thus to find the maximal modular subset of $F_n(l)$ one computes first of all $C_n(l)$. If it is non-empty, let $S_0$ be the set of all products in $F_n(l)$ such that $\overline{S}_0 = C_n(l)$. If it is not clear by other means whether $S_0$ actually contains a modular subset, then we now proceed as indicated at the beginning of this section:. Let $S_1$ be the set of all $\pi \in S_0$ such that $\pi(-1/z)$ is a linear combination of functions in $S_0$. Similarly, construct $S_2$ from $S_1$, $S_3$ from $S_2$ and so forth. Either some $S_k$ is empty, and then $F_n(l)$ contains no modular set, or $S_k = S_{k+1} \neq \emptyset$ for some $k$, and then $S_k$ is the maximal modular set in $F_n(l)$.

Assume now $n = 1$ and $l > 1$. Then $\mathbb{P}^{n-1}(\mathbb{Z}/l\mathbb{Z})$ contains only one point $[a]$. If this point yields a premodular set, one has
$$\beta_l([a]) = \mathbb{B}_2\left(\frac{1}{l}\right) = \frac{1}{6l^2}.$$

For the first equality we used that $\mathbb{B}_2(x)$ is even, decreasing between 0 and $\frac{1}{2}$, and that $\gcd(a, l) = 1$. Rewriting this identity as $5 - 6l + l^2 = 0$ we find $l = 5$ as the only solution $> 1$. And indeed, $F_1(5)$ equals $AG_5$.

Let $n = 2$ or $n = 3$. We determine, for all $l$, all non-empty $C_n(l)$. If $C_n(l)$ is non-empty, then $C_n(t)$ is non-empty for all divisors $t$ of $l$. Theorem 4.1 applied with $l$ equal to a prime $p$ shows that $C_2(p) = \emptyset$ for $p > 37$, and $C_3(p) = \emptyset$ for $p > 113$. A computer search shows that actually $C_2(p) \neq \emptyset$ only for $p = 2, 5, 7$, and $C_3(p) \neq \emptyset$ only for $p = 3, 5, 7$. Next, for each of these primes $p$, we look for powers $p^r$ such that $C_n(p^r)$ is nonempty. The possible values of $r$ are bounded by Theorem 4.1. Again by a computer search, we find that $C_2(l) \neq \emptyset$ implies $l \mid 2 \cdot 5 \cdot 7$, and that $C_3(l) \neq \emptyset$ implies $l \mid 3^2 \cdot 5 \cdot 7$. A final computer search yields then the Table 2. The above procedure to pass from premodular sets to the maximal modular one, i.e. to descend to $S_1$, $S_2$ etc., had actually only been applied twice in the course of our computations: to rule out certain functions for $n = 3$ and $l = 15$, and to prove that $W_7$ is modular.

# 6 The $[r]_l$ in terms of $l$-division values of the Weierstrass $\sigma$-function

Problems involving the action of $\Gamma$ on modular units are most conveniently studied using $l$-th division values of the Weierstrass $\sigma$-function (or Siegel units, as they are called in the literature). This relies on the following two facts: Firstly, the action $\Gamma$ on a Siegel unit is given by an explicit formula (Theorem 6.1). Secondly, if $\mathcal{S}$ denotes the group generated by the Siegel units, then $U/\mathcal{S}$ has exponent 2.

The transformation formulas are most naturally and easily derived by using the Jacobi group and considering the Weierstrass $\sigma$-function as Jacobi form. Since this approach cannot be found in the literature we present it here in form of an appendix. The resulting formulas, however, are classical and well-known.

For the complicated proof of the second fact, namely that $U/\mathcal{S}$ has exponent 2, see the book [K-L] and papers cited therein. For us, fortunately, it suffices to know the considerably simpler fact that $U/\mathcal{S}$ is a torsion group (Theorem 6.2). Since we do not know any reference to an easy and direct proof of this, we shall give such a proof here. Finally, we shall describe below the relation between Siegel units and the functions $[r]_l$ and we shall deduce from this and the two theorems on Siegel units the facts (Lemma 6.1) which were used in the preceding paragraphs without proofs.

For a row vector $\alpha = (a_1, a_2) \in \mathbb{Q}^2$, $\alpha \notin \mathbb{Z}^2$, set

$$s_\alpha = q^{-\mathbb{B}_2(a_1)/2} \prod_{\substack{n \equiv a_1 \, (\mathbb{Z}) \\ n \geq 0}} (1 - q^n e(a_2))^{-1} \prod_{\substack{n \equiv -a_1 \, (\mathbb{Z}) \\ n > 0}} (1 - q^n e(-a_2))^{-1}. \qquad (2)$$

Here $e(\dots) = \exp(2\pi i \dots)$. Note that the first product has to be taken over all non-negative $n$, whereas the second one is over strictly positive $n$ only. Moreover, $s_\alpha$ depends only on $\alpha \bmod \mathbb{Z}^2$. Finally, $s_\alpha$ has clearly no zeros or poles in the upper half plane. The functions $s_\alpha^{-1}$ are known in the literature as Siegel units.

The following theorem will be proved in the Appendix.

**Theorem 6.1.** *For all non-integral $\alpha \in \mathbb{Q}^2$ and all $A \in \mathrm{SL}(2,\mathbb{Z})$ there exist a root of unity $c(\alpha, A)$ such that*

$$s_\alpha \circ A = c(\alpha, A) \, s_{\alpha A}.$$

*(Here $aA$ means the usual action of $A$ on the row vector $a$, and $(s_\alpha \circ A)(z) = s_\alpha(Az)$.) Moreover, for a given $\alpha$, the group of $A$ such that $c(\alpha, A) = 1$ is a congruence subgroup. In particular, $s_\alpha$ is a modular unit.*

*Remark* 6.1. The numbers $c(\alpha, A)$ define obviously a cocycle of $\Gamma = \mathrm{SL}(2,\mathbb{Z})$, i.e. one always has $c(\alpha, AB) = c(\alpha, A) \cdot c(\alpha A, B)$. The actual values of $c(\alpha, A)$ will drop out automatically of the proof given below though we do not need them. In particular, $c(\alpha, A)^{12l} = 1$ for all $A$, where $l$ is the common denominator of $a_1$ and $a_2$.

**Theorem 6.2.** *Let $f$ be a modular unit on the principal congruence subgroup $\Gamma(l)$. Then a suitable positive integral power of $f$ is, up to multiplication by a constant, contained in the group generated by the Siegel units $s_\alpha$, where $\alpha$ runs through all pairs of rational numbers of the form $(\frac{r}{l}, \frac{s}{l})$ with integers $r, s$ such that $\gcd(r, s, l) = 1$.*

*Proof.* Let $U[\Gamma(l)]/\mathbb{C}^*$ be the group of modular units on $\Gamma(l)$ modulo multiplication by constants. Since the map

$$U[\Gamma(l)]/\mathbb{C}^* \to \mathbb{Z}[\Gamma(l)\backslash \mathbb{P}(\mathbb{Q})], \quad \pi \mapsto \text{divisor of } \pi$$

is injective and takes its image in the subgroup of divisors of degree 0, we conclude that the rank of $U[\Gamma(l)]/\mathbb{C}^*$ is $\leq R-1$, where $R = |\Gamma(l)\backslash\mathbb{P}(\mathbb{Q})|$ is the number of cusps of $\Gamma(l)$. It is well-known that $R$ equals the cardinality of $I$, where $I$ is a complete set of representatives for

$$\{(\frac{r}{l}, \frac{s}{l}) \in \frac{1}{l}\mathbb{Z}^2 \colon \gcd(r,s,l) = 1\}$$

modulo $\mathbb{Z}^2$ and modulo multiplication by $\pm 1$.

On the other hand, for suitable large integers $N$ the powers $s_\alpha^N$ with $\alpha \in I$ are elements in $U[\Gamma(l)]$ (in fact, one may take $N = 12l$). This is an immediate consequence of Theorem 6.1 and the remark subsequent to it.

Moreover, a relation

$$\prod_{\alpha \in I} s_\alpha^{c(\alpha)} = \text{const.}$$

holds true if and only if $c(\alpha)$, as function of $\alpha$, is constant. Indeed, a constant function $c(\alpha)$ yields a modular unit a power of which is invariant under $\mathrm{SL}(2,\mathbb{Z})$ by Theorem 6.1. Since $\mathrm{SL}(2,\mathbb{Z})$ has only one cusp, this unit must be a constant. That there is no other relation can, *e.g.*, be verified by looking at the logarithmic derivatives of the $s_\alpha$, which, by well-known theorems, span the space of Eisenstein series on $\Gamma(l)$ [H, pp. 468]. But the dimension of this space is $R-1$. Hence the rank of the subgroup of $U[\Gamma(l)]/\mathbb{C}^*$ generated by the $\mathbb{C}^* \cdot s_\alpha^N$ ($\alpha \in I$) equals $R-1$. We deduce from this that $U[\Gamma(l)]/\mathbb{C}^*$ has full rank $R-1$, and that the $\mathbb{C}^* \cdot s_\alpha^N$ ($\alpha \in I$) generate a subgroup of finite index. $\square$

**Lemma 6.1.** *Let $l \geq 1$ be an integer. For each integer $r$ not divisible by $l$ one has*

$$[r]_l = \prod_{s \bmod l} s_{(r,s)/l}.$$

*In particular, $[r]_l$ is a modular unit. For each $A = (a,b;c,d) \in \Gamma$ one has*

$$[r]_l \circ A \in c\, q^{-\frac{t^2}{2l}\mathbb{B}_2(\frac{ar}{t})}(1 + q^{1/l}\, K[\![q^{1/l}]\!]),$$

*where $t = \gcd(c,l)$, where $K$ denotes the field of $l$-th roots of unity, and where $c$ is a constant.*

*Proof.* The formula expressing $[r]_l$ in terms of the $s_\alpha$ is a simple consequence of the polynomial identity

$$\prod_{k \bmod l} (1 - e(k/l)\, Z) = 1 - Z^l.$$

By Theorem 6.1 the function $[r]_l$ is then a modular unit. The last assertion follows from the given formula, Theorem 6.1, and on using

$$-\frac{1}{2}\sum_{s \bmod l} \mathbb{B}_2\left(\frac{ra+cs}{l}\right) = -\frac{t}{2}\sum_{\substack{y \bmod l \\ y \equiv ar \bmod t}} \mathbb{B}_2\left(\frac{y}{l}\right) = -\frac{t^2}{2l}\mathbb{B}_2\left(\frac{ar}{l}\right).$$

Here the second identity is the well-known distribution property of the Bernoulli polynomial $\mathbb{B}_2(x)$. □

# 7  Appendix: The Weierstrass $\sigma$-function as Jacobi form

For $z \in \mathbb{H}$ and $x \in \mathbb{C}$ let

$$\begin{aligned}\phi(z,x) &= 2\pi i\, \eta(z)^2\, e\left(\frac{\eta'}{\eta}(z)\, x^2\right) x \prod_{\substack{l \in \mathbb{Z}z+\mathbb{Z} \\ l \neq 0}} \left(1-\frac{x}{l}\right)\exp\left(\frac{x}{l}+\frac{1}{2}\left(\frac{x}{l}\right)^2\right) \\ &= q^{\frac{1}{12}}\left(\zeta^{1/2}-\zeta^{-1/2}\right)\prod_{n \geq 1}\left(1-q^n\zeta\right)\left(1-q^n\zeta^{-1}\right),\end{aligned}$$

with $\eta(z) = q^{1/24}\prod_{n \geq 1}(1-q^n)$ denoting the Dedekind $\eta$-function, and with $\zeta^r(x) = \exp(2\pi i r x)$ (see, e.g., [S, pp. 143] for a proof of the equality of the two expressions for $\phi$). Note that $\phi(z,x)$ is, up to the factors involving $\eta$, the Weierstrass $\sigma$-function.

**Theorem 7.1.** *For* $A = \begin{pmatrix} a & b \\ c & d \end{pmatrix} \in \mathrm{SL}(2,\mathbb{Z})$ *and* $\alpha = (a_1, a_2) \in \mathbb{Z}^2$ *one has*

$$\phi\left(Az, \frac{x}{cz+d}\right) e\left(-\frac{cx^2}{2(cz+d)}\right) = \epsilon(A)\phi(z,x),$$

$$\phi(z, x+a_1 z + a_2)\, q^{a_1^2/2}\zeta^{a_1}\, (-1)^{a_1+a_2} = \phi(z,x),$$

*where* $\epsilon(A) = \eta^2(Az)/\bigl((cz+d)\eta^2(z)\bigr)$.

*Proof.* The first formula is an immediate consequence of the first definition for $\phi$ on using the identity $\mathbb{Z}\, Az + \mathbb{Z} = \frac{1}{cz+d}(\mathbb{Z}z + \mathbb{Z})$ and

$$\frac{\eta'}{\eta}(Az)\, \frac{1}{(cz+d)^2} = \frac{c}{2(cz+d)} + \frac{\eta'}{\eta}(z),$$

which in turn follows immediately from

$$\eta^2(Az) = \epsilon(A)\eta^2(z)\,(cz+d)$$

with a certain constant $\epsilon(A)$. The second transformation formula can be directly checked using the second formula for $\phi$. □

It is convenient to interpret these transformation laws as an invariance property with respect to a certain group, namely the Jacobi group $\mathbb{J}(\mathbb{Z})$. For a ring $R$ (commutative, with 1) denote by $\mathbb{J}(R)$ the group of all triples $(A, \alpha, n)$ of matrices $A \in \mathrm{SL}(2, R)$, row vectors $\alpha \in R^2$ and $n \in R$, equipped with the multiplication law

$$(A, \alpha, n) \cdot (B, \beta, n') = (AB, \alpha B + \beta, n + n' + \det\begin{pmatrix}\alpha B\\ \beta\end{pmatrix}).$$

The Jacobi group $\mathbb{J}(\mathbb{R})$ acts on functions $\psi(z,x)$ defined on $\mathbb{H} \times \mathbb{C}$ by

$$(\psi|(A,\alpha,n))(z,x)$$
$$= e_2\left(-\frac{cx^2}{(cz+d)} + a_1^2 z + 2a_1 x + a_1 a_2 + n\right) \psi\left(Az, \frac{x + a_1 z + a_2}{cz+d}\right)$$

(with $A$ and $\alpha$ as in the above Lemma, and with $e_2(\dots) = e(\frac{1}{2}[\dots])$). That this is indeed an action can be verified by a direct (though subtle) computation [E-Z, Theorem 1.4]. Using this group action the formulas of Lemma can now be reinterpreted as

$$\phi|g = \rho(g)\phi,$$

for all $g \in \mathbb{J}(\mathbb{Z})$, where

$$\rho((A,\alpha,n)) = (-1)^{a_1+a_2+a_1 a_2+n}\epsilon(A).$$

From this transformation law for $\phi$ it is clear that $\rho$ defines a character of $\mathbb{J}(\mathbb{Z})$, as can, of course, also be checked directly.

*Proof of Theorem 6.1.* For $\alpha = (a_1, a_2) \in \mathbb{Q}^2$ and $\beta = (b_1, b_2) \in \mathbb{Z}^2$ we have

$$\phi|(1, \beta+\alpha, 0) = \phi|[(1,\beta,0) \cdot (1,\alpha,0) \cdot (1, 0, \det\begin{pmatrix}\alpha\\ \beta\end{pmatrix})]$$
$$= \rho((1,\beta,0))\, e_2\left(\det\begin{pmatrix}\alpha\\ \beta\end{pmatrix}\right)\phi|(1,\alpha,0).$$

Call the factor in front of $\phi|(1,\alpha,0)$ on the right hand side $C(\alpha,\beta)$. It can easily be checked that
$$C(\alpha,\beta) = \delta(\alpha+\beta)/\delta(\alpha),$$
where
$$\delta(\alpha) = -\rho((1,\lfloor\alpha\rfloor,0))\, e_2\left(\det\begin{pmatrix}\alpha\\ \lfloor\alpha\rfloor\end{pmatrix}\right)\, e_2\left(-(a_2-\lfloor a_2\rfloor)(a_1-\lfloor a_1\rfloor-1)\right),$$
with $\lfloor\alpha\rfloor = (\lfloor a_1\rfloor, \lfloor a_2\rfloor)$. Thus, if we set, for $\alpha \in \mathbb{Q}^2$, $\alpha \notin \mathbb{Z}^2$,
$$S_\alpha = \delta(\alpha)\phi^{-1}|(1,\alpha,0),$$
then $S_\alpha = S_{\alpha+\beta}$ for $\beta \in \mathbb{Z}^2$. From the transformation law for $\phi$ under $\mathbb{J}(\mathbb{Z})$ we obtain
$$S_\alpha|(A,0,0) = \epsilon(A)^{-1}\frac{\delta(\alpha)}{\delta(\alpha A)}S_{\alpha A}.$$
A simple calculation shows that $s_\alpha(z) = S_\alpha(z,0)$. From this Theorem 6.1 follows. $\square$

## Acknowledgments.

We would like to thank Don Zagier and Herbert Gangl for discussions and comments, and Michael Rösgen for comments on an early draft version of this article. Nils-Peter Skoruppa thanks the Department for Applied Mathematics and Theoretical Physics for financial support during a stay in Cambridge where parts of this article were written. Wolfgang Eholzer was supported by the EPSRC; partial support from PPARC and from EPSRC, grant GR/J73322, is also acknowledged.

Wolfgang Eholzer  
Deutsche Börse Systems AG  
Neue Börsenstrasse 1  
60487 Frankfurt, Germany  
wolfgang.eholzer@deutsche-boerse.com

Nils-Peter Skoruppa  
Universität Siegen, Fachbereich Mathematik  
Walter-Flex-Strasse 3  
57068 Siegen, Germany  
skoruppa@math.uni-siegen.de